\documentclass[11pt,reqno]{amsart}
\usepackage{amssymb,amsmath}\newcommand{\be}{\begin{eqnarray}}
\newcommand{\ee}{\end{eqnarray}}

\newcommand{\R}{{\mathbb R}}

\newcommand{\E}{{\bf E\,}}

\newtheorem{theorem}{Theorem}

\theoremstyle{definition}
\theoremstyle{remark}\numberwithin{equation}{section}\input epsf.sty
\begin{document}\thispagestyle{empty}

\title[Orthogonal martingales]{Subordination by orthogonal  martingales in $L^{p}, 1<p\le 2$}
\author{Prabhu Janakiraman}\address{Prabhu Janakiraman, Dept. Math., Purdue Univ.,
{\tt pjanaki@math.purdue.edu}}
\author{Alexander Volberg}\address{Alexander Volberg, Dept. of  Math., Mich. State Univ.,
{\tt volberg@math.msu.edu}}


\maketitle

\section{Introduction: Orthogonal martingales and the Beurling-Ahlfors transform}

We are given two martingales on the filtration of the two dimensional Brownian motion. One is subordinated to another. We want to give an estimate of $L^p$-norm of a subordinated one via the same norm of a dominating one. In this setting this was done by Burkholder in \cite{Bu1}--\cite{Bu8}. If one of the martingales is orthogonal, the constant should drop. This was demonstrated in \cite{BaJ1}, when the orthogonality is attached to the subordinated martingale and when $2\le p<\infty$. This note contains an (almost obvious) observation that the same idea can be used in the case when the orthogonality is attached to a dominating martingale and $1<p\le 2$. Two other complementary regimes are considered in \cite{BJV_La}. When both martingales are orthogonal, see \cite{BJV_Le}. In these two papers the constants are sharp.
We are not sure of the sharpness of the constant in the present note.

A complex-valued martingale $Y=Y_1+iY_2$ is said to be $orthogonal$ if the quadratic variations of the coordinate martingales are equal and their mutual covariation is $0$:
\[\left<Y_1\right>=\left<Y_2\right>,\hspace{3mm} \left<Y_1,Y_2\right>=0.\]
In \cite{BaJ1}, Ba\~nuelos and Janakiraman make the observation that the martingale associated with the Beurling-Ahlfors transform is in fact an orthogonal martingale. They show that Burkholder's proof in \cite{Bu3} naturally accommodates for this property and leads to an improvement in the estimate of $\|B\|_p$.
\begin{theorem}\label{BurkBaJa} (One-sided orthogonality as allowed in Burkholder's proof)
\begin{enumerate}
\item (Left-side orthogonality) Suppose $2\leq p<\infty$. If $Y$ is an orthogonal martingale and $X$ is any martingale such that $\left<Y\right>\leq \left<X\right>$, then 
\begin{equation}\label{estBJ}
\|Y\|_p \leq \sqrt{\frac{p^2-p}{2}}\|X\|_p.
\end{equation}
\item (Right-side orthogonality) Suppose $1<p<2$. If $X$ is an orthogonal martingale and $Y$ is any martingale such that $\left<Y\right>\leq \left<X\right>$, then 
\begin{equation}\label{estBJ1}
\|Y\|_p \leq \sqrt{\frac{2}{p^2-p}}\|X\|_p.
\end{equation}
\end{enumerate}It is not known whether these estimates are the best possible.
\end{theorem}

\noindent{\bf Remark.}
The result for right-side orthogonality is stated in \cite{JVV} and not in \cite{BaJ1}. But \cite{JVV} has a complicated (though funny and interesting) proof by construction a family of new Bellman functions very different from the original Burkholder's function.    The goal of this small note is to demonstrate how one adapt the idea of \cite{BaJ1} to the right-orthogonality and $1<p\le 2$ regime. We use just a well-known Burkholder's function here, exactly along the lines of \cite{BaJ1}.

If $X$ and $Y$ are the martingales associated with $f$ and $Bf$ respectively, then $Y$ is orthogonal, $\left<Y\right>\leq 4\left<X\right>$ and hence by (1), we obtain
\begin{equation}\label{Beurest}
\|Bf\|_p\leq \sqrt{2(p^2-p)}\|f\|_p \textrm{ for } p\geq 2.
\end{equation}
By interpolating this estimate $\sqrt{2(p^2-p)}$ with the known $\|B\|_2=1$, Ba\~nuelos and Janakiraman establish the present best estimate in publication:
\begin{equation}\label{best-est}
\|B\|_p\leq 1.575 (p^*-1).
\end{equation}

\section{New Questions and  Results}
Since $B$ is associated with left-side orthogonality and since we know $\|B\|_p = \|B\|_{p'}$, two important questions are
\begin{enumerate}
\item If $2\leq p<\infty$, what is the best constant $C_p$ in the left-side orthogonality problem: $\|Y\|_p\leq C_p\|X\|_p$, where $Y$ is orthogonal and $\left<Y\right>\leq \left<X\right>$?
\item Similarly, if $1<p'<2$, what is the best constant $C_{p'}$ in the left-side orthogonality problem?
\end{enumerate}
We have separated the two questions since Burkholder's proof (and his function) already gives a good answer when $p\geq 2$. It may be (although we have now some doubts about that) the best possible as well. However no estimate (better than $p-1$) follows from analyzing Burkholder's function when $1<p'<2$. Perhaps, we may hope, $C_{p'}<\sqrt{\frac{p^2-p}{2}}$ when $1<p'=\frac{p}{p-1}<2$, which would then imply a better estimate for $\|B\|_p$. This paper 'answers' this hope in the negative by finding $C_{p'}$; see Theorem \ref{MAINThm}. We also ask and answer the analogous question of right-side orthogonality when $2<p<\infty$. In the spirit of Burkholder \cite{Bu8}, we believe these questions are of independent interest in martingale theory and may have deeper connections with other areas of mathematics.

\bigskip 

\noindent{\bf Remark.}
The following sharp estimates are proved in \cite{BJV_La}, they cover the left-side orthogonality for the regime $1<p\le 2$ and the right-side orthogonality for the regime $2\le p<\infty$. Notice that two complementary regimes have the estimates: for $2\le p<\infty$ and left-side orthogonality in \cite{BaJ1}, for $1<p\le 2$ in this note and in \cite{JVV}, but the sharpness is dubious.

\bigskip

\begin{theorem}\label{MAINThm}Let $Y=(Y_1,Y_2)$ be an orthogonal martingale and $X=(X_1,X_2)$ be an arbitrary martingale.
\begin{enumerate}
\item Let $1<p' \leq 2$. Suppose $\left<Y\right>\leq\left<X\right>$. Then the least constant that always works in the inequality $\|Y\|_{p'}\leq C_{p'}\|X\|_{p'}$ is 
\begin{equation}\label{constp<2}
C_{p'} = \frac{1}{\sqrt{2}}\frac{z_{p'}}{1-z_{p'}}
\end{equation}
where $z_{p'}$ is the least positive root in $(0,1)$ of the bounded Laguerre function $L_{p'}$.
\item Let $2\leq p<\infty$. Suppose $\left<X\right>\leq \left<Y\right>$. Then the least constant that always works in the inequality $\|X\|_{p}\leq C_{p}\|Y\|_{p}$ is 
\begin{equation}\label{constp>2}
C_{p} = \sqrt{2}\frac{1-z_{p}}{z_{p}}
\end{equation}
where $z_{p}$ is the least positive root in $(0,1)$ of the bounded Laguerre function $L_{p}$.
\end{enumerate}
\end{theorem}
The Laguerre function $L_p$ solves the ODE 
\[sL_p''(s)+(1-s)L_p'(s)+pL_p(s) = 0.\] These functions are discussed further and their properties deduced in section \eqref{laguerresection}; see also \cite{BJV}, \cite{C}, \cite{CL}.

As mentioned earlier, (based however on numerical evidence) we believe in general $\sqrt{\frac{p^2-p}{2}}<C_{p'}<p-1$ and that these theorems cannot imply better estimates for $\|B\|_p$. However based again on numerical evidence, the following conjecture is made.

\noindent{\bf Conjecture.}
For $1<p'=\frac{p}{p-1}<2$, $C_{p'} = C_p$, or equivalently,
\[ \frac{1}{\sqrt{2}}\frac{z_{p'}}{1-z_{p'}} = \sqrt{2}\frac{1-z_{p}}{z_{p}}.\]

It is conjecture relating the roots of the Laguerre functions. Notice that such a statement is not true with the constants from Theorem \ref{BurkBaJa}, and $\sqrt{\frac{2}{p'^2-p'}}<\sqrt{\frac{p^2-p}{2}}$ for all $p>2$. So this conjecture (if true) suggests some distinct implications for the two settings. Note on the other hand, that the form of the two sets of constants are very analogous.

\section{Right-side orthogonality, $1<p\le 2$ regime, Burkholder's function}

We just repeat the approach of \cite{BaJ1}. Let 
$$
\alpha_p := p\bigg(1-\frac1{p^*}\bigg)^{p-1}\,, 1<p\le 2\,.
$$
For $x\in R^2, y\in R^2$ we define following Burkholder:
$$
v(x,y) := \|y\|^p-(p^*-1)^{p}\|x\|^p\,.
$$
We consider Burkholder's function
$$
u(x,y) := \alpha_p (\|y\|-(p^*-1)\|x\|)(\|x\|+\|y\|)^{p-1}\,.
$$
Then ($1<p\le 2$)
$$
(p-1) u(x,y) =-\alpha_p(\|x\|-(p-1)\|y\|)(\|x\|+\|y\|)^{p-1}\,.
$$
So if we denote $G(t) := u(x+ht, y+kt)$ we have
$$
G''(0) = -\alpha_p (A+B+C)\,,
$$
where 
$$
A= p(p-1)(\|h\|^2-\|k\|^2) (\|x\|+\|y\|)^{p-1}\,,
$$
$$
B= (2-p)p(\|h\|^2-(\frac{x}{\|x\|},h)^2)\|x\|^{-1}(\|x\|+\|y\|)^{p-1}\,.
$$
And $C\ge 0$.

Also $(p-1)u(x,y)\le 0$ if $\|y\|\le \|x\|$.

Now let temporarily  $X_t=(X^1_t, X^2_t), Y_t=(Y^1_t, Y^2_t)$ denote two $\R^2$--valued martingales on the filtration of $2$--Brownian motion, and let
\begin{equation}
\label{locort3}
d\langle X^1, X^2\rangle = h^1_1h^2_1 + h^1_2h^2_2 =0\,.
\end{equation}
\begin{equation}
\label{equalnorms4}
d\langle X^1, X^1\rangle = (h^1_1)^2 + (h^1_2)^2=
d\langle X^2, X^2\rangle=(h^2_1)^2 + (h^2_2)^2\,.
\end{equation}
And let us have the following subordination by the orthogonal martingale assumption:
\begin{equation}
\label{subord5}
d\langle Y,Y\rangle \le \frac{p}{2(p-1)}d\langle X, X\rangle\,,
\end{equation}
or
\begin{equation}
\label{subord6}
(k^1_1)^2 + (k^1_2)^2+(k^2_1)^2 + (k^2_2)^2\le \frac{p}{2(p-1)}((h^1_1)^2 + (h^1_2)^2+(h^2_1)^2 + (h^2_2)^2)\,,
\end{equation}

We write It\^o's formula for $\E u(X_t, Y_t)$:

$$
\E u(X_t, Y_t) =\E u(X_0, Y_0) -\frac{\alpha_p}{2}\E\int_0^t (A(t) +B(t) +C(t))\,dt\,,
$$ where (see above) 
$$
A(t)= p(p-1)(d\langle X, X\rangle_t-d\langle Y, Y\rangle_t) (\|X_t\|+\|Y_t\|)^{p-1}\,,
$$
$$
B= (2-p)p(d\langle X, X\rangle_t-[(\frac{X_t}{\|X_t\|},\overrightarrow{H_1})^2)+(\frac{X_t}{\|X_t\|},\overrightarrow{H_2})^2)]\|X_t\|^{-1}(\|X_t\|+\|Y_t\|)^{p-1}\,.
$$
And $C(t)\ge 0$.

Here we denote
$$
H_1 = (h^1_1, h^2_1)\, , H_2 =(h^2_1, h^2_2)\,,
$$
or we can say that $H_1$ is a ``vector of $x$ stochastic derivatives of vector process $X$" and $H_2$ is a ``vector of $y$ stochastic derivatives of vector process $X$".
By \eqref{locort3} and \eqref{equalnorms4} we get that the expression in $[\cdot]$ is
$$
[(\frac{X_t}{\|X_t\|},\overrightarrow{H_1})^2)+(\frac{X_t}{\|X_t\|},\overrightarrow{H_2})^2)] =\frac1{2}d\langle X,X\rangle\,.
$$
Hence, if $\|Y_0\|\le \|X_0\|$ we get (as $\|X_t\|^{-1}(\|X_t\|+\|Y_t\|)^{p-1} \ge (\|X_t\|+\|Y_t\|)^{p-2}$)
$$
\E u(X_t, Y_t) \le -\frac{\alpha_p}{2}\E\int_0^t \{p(p-1)(d\langle X,X\rangle -d\langle Y,Y\rangle +\frac{2-p}{2(p-1)} d\langle X,X\rangle)\}\,dt\,,
$$
or 
$$
\E u(X_t, Y_t) \le -\frac{\alpha_p}{2}\E\int_0^t \{p(p-1)\bigg(\frac{p}{2(p-1)}d\langle X,X\rangle -d\langle Y,Y\rangle \bigg)\}\,dt\le 0\,,
$$
because the integrand is positive: see the assumption of subordination \eqref{subord5}. Therefore, using Burkholdr's discovery that
 $$
 v(x,y) = \|y\|^p-(p^*-1)^{p}\|x\|^p\le u(x,y)
 $$
 we get
 $$
 \E( \|Y_t\|^p-(p^*-1)^{p}\|X_t\|^p\le\E u(X_t,Y_t)\le 0\,,
 $$ and we obtain
 $$
 \|Y\|_p \le (p^*-1)\|X\|_p\,.
 $$ 
 Consider $\widetilde{X} := \sqrt{\frac{p}{2(p-1)}} X$. Then \eqref{locort3} means ortogonality 
 $d\langle \widetilde{X}^1, \widetilde{X}^2\rangle=0$. Assumption \eqref{equalnorms4} means $d\langle\widetilde{X}^1, \widetilde{X}^1\rangle = d\langle \widetilde{X}^2, \widetilde{X}^2\rangle$, and \eqref{subord5} means
 $d\langle Y, Y\rangle \le d\langle \widetilde{X}, \widetilde{X}\rangle$. Changing $X$ to $\widetilde{X}$ we see that we proved
 \begin{theorem}
 \label{1p2th}
 Let $1<p\le 2$, let $X_t, Y_t$ be two martingales on the filtration of $2$--dimensional Brownian motion.
 Let $X$ be an orthogonal martingale, namely $d\langle X^1, X^2\rangle=0$ and $d\langle X^1, X^1\rangle=d\langle X^2, X^2\rangle$. Suppose that $Y$ is subordinated to $X$:
 $$
 d\langle Y, Y\rangle\le d\langle X, X\rangle\,.
 $$
 Then
 $$
 \|Y\|_p \le \sqrt{\frac{2}{p^2-p}}\|X\|_p\,.
 $$
 \end{theorem}

\markboth{}{\sc \hfill \underline{References}\qquad}


\begin{thebibliography}{XXXXXX}
\label{rf}



\bibitem[BaJ1]{BaJ1}{\sc R. Banuelos, P. Janakiraman},
{\em $L^p$--bounds for the  Beurling--Ahlfors transform}. Trans.
Amer. Math. Soc. {\bf 360} (2008), no. 7, 3603--3612.

\bibitem[Bu1]{Bu1}
 {\sc D.~Burkholder}, {\em Boundary value problems and sharp estimates for the martingale transforms}, Ann. of Prob. {\bf 12} (1984), 647--702.
 
 \bibitem[Bu2]{Bu2}
 {\sc D.~Burkholder}, {\em An extension of classical martingale inequality}, Probability Theory and Harmonic Analysis, ed. by J.-A. Chao and W. A. Woyczynski, Marcel Dekker, 1986.

\bibitem[Bu3]{Bu3}
 {\sc D.~Burkholder}, {\em Sharp inequalities for martingales and stochastic
integrals}, Colloque Paul L\' evy sur les Processus Stochastiques
(Palaiseau, 1987), Ast\' erisque No. 157-158 (1988), 75--94.

\bibitem[Bu4]{Bu4}
 {\sc D.~Burkholder}, {\em Differential subordination of harmonic functions and martingales}, (El Escorial 1987)
 Lecture Notes in Math., {\bf 1384} (1989), 1--23.
 
 \bibitem[Bu5]{Bu5}
 {\sc D.~Burkholder}, {\em Explorations of martingale theory and its applications}, Lecture Notes in Math. {\bf 1464} (1991), 1--66.
 
 
 \bibitem[Bu6]{Bu6}
 {\sc D.~Burkholder}, {\em Strong differential subordination and stochastic integration}, Ann. of Prob. {\bf 22} (1994), 995--1025.

\bibitem[Bu7]{Bu7}
 {\sc D.~Burkholder}, {\em A proof of the Peczynski's conjecture for the Haar
system}, Studia Math., {\bf 91} (1988), 79--83.

\bibitem[Bu8]{Bu8}
{\sc D.~Burkholder}, {\em Martingales and Singular Integrals in Banach Spaces}, Handbook of the Geometry of Banach Spaces, Vol. 1, Chp. 6., (2001), 233-269.



\bibitem[BJVLa]{BJV_La}  {\sc A. Borichev, P. Janakiraman,  A. Volberg}, {\em Subordination by orthogonal  martingales in $L^{p}$ and zeros of Laguerre polynomials}
Preprint, 2009, pp. 1--27, sashavolberg.wordpress.com


\bibitem[BJVLe]{BJV_Le}  {\sc A. Borichev, P. Janakiraman,  A. Volberg}, {\em On Burkholder function for orthogonal  martingales  and zeros of Legendre polynomials}, arXiv:1002.2314,
Preprint, 2009, pp. 1--36, sashavolberg.wordpress.com



\bibitem[JVV]{JVV} {\sc P. Janakiraman, V. Vasyunin, A. Volberg}, {\em Some new Bellman functions and subordination by orthogonal martingales in $L^p$, $1< p\le 2$}, Preprint, 2009, pp. 1--30, sashavolberg.wordpress.com

\end{thebibliography}
\end{document}